\documentclass[11pt]{elsart}

\usepackage{eucal,amssymb,enumerate}
\usepackage{graphics}

\journal{Journal of Algebra}

\input xy
\xyoption{all}

\newtheorem{thrm}{Theorem}[section]

\newtheorem{lemma}[thrm]{Lemma}

\newtheorem{prp}[thrm]{Proposition}

\newtheorem{crl}[thrm]{Corollary}

\newtheorem{dfn}[thrm]{Definition}

\newtheorem{exm}[thrm]{Example}

\begin{document}

\begin{frontmatter}
\title{Structure on the set of closure operations of a commutative ring}
\author{Janet C. Vassilev}
\address{The University of California, Riverside, CA 92521}
\ead{jvassil@ucr.edu} \ead[url]{http:\slash\slash
math.ucr.edu\slash$^\sim$jvassil}

\begin{abstract}{We investigate the algebraic structure on the set
of closure operations of a ring.
We show the set of closure operations is not a monoid under
composition for a discrete valuation ring.  Even the set of
semiprime operations over a DVR is not a monoid; however, it is the
union of two monoids, one being the left but not right act of the
other. We also determine all semiprime operations over the ring
$K[[t^2, t^3]]$.
 }
 \end{abstract}

\end{frontmatter}

\section{Introduction}

Let $I \mapsto I_c$ be an operation on the set of ideals of a ring $R$.
Consider the following properties where $I$ and $J$ are ideals and $b$ is
a regular element:
\begin{enumerate}[(a)]
\item $I \subseteq I_c$
\item If $I \subseteq J$, then $I_c \subseteq J_c$.
\item $(I_c)_c=I_c$.
\item $I_cJ_c \subseteq (IJ)_c$
\item $(bI)_c=bI_c$
\end{enumerate}

If $I \rightarrow I_c$ satisfies (a)-(c) above, we call $I \rightarrow
I_c$ a {\it closure operation}.  If $I \rightarrow I_c$ is a closure
operation and also satisfies (d) above, we call $I \rightarrow I_c$ a {\it
semiprime operation}.  If $I \rightarrow I_c$ is semiprime and also
satisfies (e), then we say $I \rightarrow I_c$ is a {\it prime operation}.

The definition of prime operation or $\prime$-operation for the set
of fractional ideals of an integral domain was given by Krull in his
1935 book, Idealtheorie \cite{Kr1}.  In his original definition, he
actually added a sixth property (f) $I_c+J_c \subseteq (I+J)_c$.
Then in his 1936 paper \cite{Kr2}, he discusses the integral
completion or $b$-operation in terms of $\prime$-operations on the
set of fractional ideals and mentions that he left out the
properties (g) $R=R_c$ and (h) $(I_c \cap J_c)_c=I_c \cap J_c$.  In
fact, Sakuma \cite{Sa} shows in 1957 that when looking at prime
operations on the set of fractional ideals of a domain, properties
(d), (f) and (h) are consequences of properties of (a), (b), (c),
(e) and (g). In 1964, Petro \cite{Pe} called the operations
satisfying properties (a)-(d) on the set of fractional ideals
semiprime operations.   The first reference to integral closure
strictly on the set of ideals of a commutative ring seems to be
Northcott and Rees' 1954 paper on reductions \cite{NR}. In 1969,
Kirby \cite{Ki} seems to be the first to discuss general closure
operations on the set of ideals over a commutative ring with
identity. The terms prime and semiprime operation were reintroduced
on the set of ideals of a commutative ring by Ratliff in his 1989
paper \cite{Ra} on $\Delta$-closures of ideals. Heinzer, Ratliff and
Rush \cite{HRR} also use the term semiprime operation when referring
to the basically full closure on the set of $m$-primary submodules
of a module over a local ring $(R,m)$.

There are many well known closure operations defined on a commutative
ring, such as:  integral closure, tight closure if the ring contains a
field \cite{HH}, $\Delta$-closure \cite{Ra}, basically full closure
\cite{HRR}, etc.  It is known that all of these closure operations are
contained in the integral closure, excluding the $\Delta$-closure.
However, if $\Delta$ doesn't contain any ideals which are contained in a
minimal prime, then the $\Delta$-closure is also contained in the integral
closure. Otherwise, the relationships between these other closures is not
as well understood.  Knowing the structure on the set of closure
operations may shed some light on this relationship.

Abstractly, closure operations are elements of the monoid of maps
from the set of ideals, $\mathfrak{I}$, of a ring to itself,
$M_{\frak{I}}=\{f:\mathfrak{I} \rightarrow \mathfrak{I}\}$
satisfying the above properties. For example, $C_R$ is the set of
maps satisfying (a)-(c), $S_R$ is the set of maps satisfying (a)-(d)
and $P_R$ is the set of maps satisfying (a)-(e). $C_R$, $S_R$ and
$P_R$ are all partially ordered sets, but otherwise these sets are
in general poorly behaved. In Section 2, we will give examples
showing that $C_R$ is not even a monoid in the nice case that $R$ is
a discrete valuation ring. Then in Section 3 we show that $S_R$ for
a discrete valuation ring $R$ is almost a monoid. In fact, $S_R$ is
the union of two submonoids of $M_{\frak{I}}$, one a left but not a
right act of the other. Also we show that $P_R$ is a monoid. We are
also able to extend our results to semiprime and prime operations
over a Dedekind domain. In Section 4, we consider closure operations
over the semigroup ring $K[[t^2,t^3]]$ and determine all the
semiprime operations over $K[[t^2,t^3]]$.

I would like to thank Dave Rush, for giving some inspirational talks on
closure operations in the commutative algebra seminar at UC Riverside and
for his many suggestions for strengthening this work. I would also like to
thank Bahman Engheta and Jooyoun Hong for being a sounding board for my
ideas and Hwa Young Lee for pointing out an error in the work over
$K[[t^2,t^3]]$.  I would also like to thank the referee for his/her many
helpful suggestions which greatly improved the reading of the paper.

\section{Preliminaries}

Recall that $(S,\circ)$ is a {\it semigroup} if $\circ$ is an
associative binary operation on $S$.  We say that a semigroup
$(S,\circ)$ is a {\it monoid} if there is a unique identity element
$e$ in $S$ such that $es=se=s$ for all $s \in S$.  In particular,
the whole numbers $\mathbb{N}_0=\{0,1,2, \ldots, n, \ldots\}$ is a
monoid under addition, with identity $0$.

Let $R$ be a commutative ring, $\frak{I}=\{I \subseteq R|I \textrm{
an ideal of } R\}$ and $M_\frak{I}=\{f:\frak{I} \rightarrow
\frak{I}\}.$ $M_\frak{I}$ is clearly a monoid under composition of
maps, with identity the identity map $e:\frak{I} \rightarrow
\frak{I}$, and function composition is associative. $C_R$ will be
the subset of $M_\frak{I}$ consisting of closure operations.  Hence
the $f_c$ in $C_R$ are the set of maps satisfying the following
three properties: (a) $f_c(I) \supseteq I$, (b) $f_c$ preserves
inclusions in $R$, and (c) $f_c\circ f_c=f_c$. $S_R$ will be the set
of semiprime operations of $R$, i.e. $S_R$ are the maps in $C_R$
which also satisfy $f_c(I)f_c(J) \subseteq f_c(IJ)$. $P_R$ will be
the set of prime operations of $R$, i.e. maps in $S_R$ which also
satisfy (e) $f_c(bI)=bf_c(I)$.  We note that if $C_R$, $S_R$ or
$P_R$ are monoids, by property (c), they will be band monoids.

\begin{dfn}  A monoid is a band monoid if every element is
idempotent. \end{dfn}

We will say $f_{c_1} \leq f_{c_2}$ for two different closure operations if
$f_{c_1}(I) \subseteq f_{c_2}(I)$ for all $I \in \frak{I}$.

\begin{prp} $C_R$,$S_R$ and $P_R$ are partially ordered sets. \end{prp}

The proof is straightforward as the ideals of $R$ are are partially
ordered under containment.

Now let us consider the algebraic structure of $C_{R}$, $R$ a
commutative ring. Unfortunately, $C_{R}$ is not a submonoid under
composition even for a discrete valuation ring.

\begin{exm}  \label{closex}   $C_R$, where $(R,P)$ is a discrete valuation ring, is not a monoid.
The ideals of $R$ have the form $P^{i}$ for all $i \geq 0$ and
$(0)$.  Let $f_n: {\frak{I}} \rightarrow {\frak{I}}$ and $g_n:
{\frak{I}} \rightarrow {\frak{I}}$ be defined as follows
 $$f_n(P^i)=\left\{ \begin{array}{l}P^i \textrm{ for } i \leq n\\
                                   P^n \textrm{ for } i
                                   >n.\end{array}\right. \textrm{ and } \quad g_n(P^i)=
                                   \left\{ \begin{array}{l} R \textrm{ for } i \leq n\\
                                   P^n \textrm{ for } i
                                   >n.\end{array}\right.$$ and
$f_n(0)=(0)=g_n(0)$.
If $m > n$, then $$f_n\circ g_m (P^i)=\left\{ \begin{array}{l}R \textrm{ for } i <m\\
                                   P^n \textrm{ for } i > m.\end{array}\right.$$
This fails property (c) as $(f_n\circ g_m)\circ(f_n\circ g_m)(P^i)=R$ for
all $i$.
\end{exm}

We will see in the next section that $g_n$ in the above example is
not a semiprime operation, because semiprime operations are not
allowed to have any finite jumps.

In Example \ref{closex} we see that the maps $f_n$ and $g_n$ are
bounded maps on the ideals of $R$. This prompts the following
definition for closure operations of commutative rings:

\begin{dfn}  We say a closure operation $f_c$ is bounded on a commutative ring $R$ if for
every maximal ideal $m$ of $R$, there is an $m$-primary ideal $I$
such that for all $m$-primary $J \subseteq I$, $f_c(J)=I$.  If this
is not the case, we will say that $f_c$ is an unbounded closure
operation.
\end{dfn}

We define bounded in this way for $m$-primary ideals, because it would be
hard to come up with a precise statement for all ideals.

\section{Algebraic structure on $S_R$ and $P_R$ when $R$ is a Dedekind domain}

It seems unlikely that $S_R$ and $P_R$ are submonoids of
$M_{\frak{I}}$ for a general commutative ring $R$, but in the case
that $R$ is a discrete valuation ring, $P_R$ is the trivial
submonoid of $M_\frak{I}$ and $S_R$ decomposes into the union of two
submonoids whose only common element is the identity. We use the
following definition to explain their relationship.

\begin{dfn} Let $S$ be a monoid and $A$ any set, then we say $A$
is a left (right) $S$-act if there is a map $\delta: S \times
A\rightarrow A$ ($\delta: A \times S\rightarrow A$) satisfying
$\delta(st,a)=\delta(s,\delta(t,a))$
($\delta(a,st)=\delta(\delta(a,s),t)$) for every $a \in A$ and $s,t
\in S$ and $\delta(e,a)=a$ ($\delta(a,e)=a$) for all $a \in A$ where
$e$ is the identity of $S$.
\end{dfn}

\begin{prp} \label{dvrsr} When $(R,P)$ is a discrete valuation ring, $S_R$ can be decomposed into
the union of two submonoids
$$M_0=\{ e\} \bigcup \left\{ f_m \in M_\frak{I} \left| f_m(P^i)=
\left\{ \begin{array}{l}P^{i} \textrm{ for } 0 \leq i <m\\
             P^m \textrm{ for } i \geq m \end{array}\right.
\textrm{ and } f_m(0)=(0) \right\}\right.  $$

and

$$M_f=\{ e\} \bigcup \left\{ g_m \in M_\frak{I} \left| g_m(P^i)=
\left\{ \begin{array}{l}P^{i} \textrm{ for } 0 \leq i <m\\
             P^m \textrm{ for } i \geq m \end{array}\right.
\textrm{ and } g_m(0)=P^m  \right\}\right.  $$

where $M_f$ is a left $M_0$-act but not a right $M_0$-act under
composition.\end{prp}

Before proving the proposition, we need the following lemma:

\begin{lemma} \label{spnc}   Let  $f_c$ be a semiprime operation
on the discrete valuation ring $(R,P)$. Then if $f_c$ is constant for
$P^i$ on a finite interval $m \leq i \leq n$ for $m <n$, then there exists
a $j \leq m$ such that $f_c(P^i)=P^j$ for all $i \geq j$.
\end{lemma}

{\bf Proof:} The ideals of $R$ have the form
 $P^i$ and they are totally ordered.  Being a closure operation, $f_c(P^i) =P^j \supseteq
P^i$, where $j \leq i$, since $f_c$ must be increasing on the ideals of
$R$.

Suppose $f_c$ is constant for $P^i$, where $m \leq i \leq n$, $m <n$. For
all such $i$ suppose that $f_c(P^i)=P^j$. Then $f_c(P^j) =
f_c(f_c(P^i))=f_c(P^i)=P^j$. Thus $P^m \subseteq f_c(P^m)=P^j$ and $j \leq
m$.  Since $f_c$ is increasing we see that $f_c(P^i)=P^j$ for all $j\leq
i<m$.

We know $f_c(P^n)=P^j$ by assumption. If we show that $f_c(P^{n+1})=P^j$,
then by induction, $f_c(P^i)=P^j$ for all $i \geq j$. Then once again, the
fact that $f_c$ is increasing implies that $f_c(P^{n+1})=P^{k}\subseteq
f_c(P^j)=P^j$ for $j \leq k \leq n+1$. Since $f_c$ is a closure operation,
$f_c(f_c(P^{n+1}))=f_c(P^k)=P^k \subseteq f_c(P^j)=P^j$. So either
$f_c(P^{n+1})=P^j$ or $f_c(P^{n+1})=P^{n+1}$. Suppose the latter.  Since,
$f_c$ is a semiprime operation, then $f_c(P^i)f_c(P^k) \subseteq
f_c(P^{i+k})$ for all $i$ and $k$; however, $f_c(P)f_c(P^{n})\supseteq
P^{j+1}$ properly contains $f_c(P^{n+1})=P^{n+1}$. Thus
$f_c(P^{n+1})=P^j$.   \qed

{\bf Proof of \ref{dvrsr}:} The ideals of a discrete valuation ring
$(R,P)$ are either of the form $P^i$ for $i \geq 0$ or $(0)$ and
they are totally ordered $R \supseteq P \supseteq P^2 \supseteq
\cdots P^m \supseteq \cdots \supseteq (0)$.

By Lemma \ref{spnc}, we know that any semiprime operation $f_c$ on $R$
which is constant on some finite interval has the property that
$f_c(P^i)=P^m$ for all $i \geq m$ for some $m$. I claim that for $i \leq
m$, $f_c(P^i)=P^i$.  Suppose not, then $f_c(P^i)=P^k$ some $k \leq i$
since $f_c$ is increasing.  Then for $k \leq j \leq i$,
$P^k=f_c(P^k)\subseteq f_c(P^j) \subseteq f_c(P^i)=P^k$. If $k < i$ then
by Lemma \ref{spnc} $f_c(P^i)=P^k$ on interval $i \geq k$ contradicting
the fact that for $i \geq m$, $f_c(P^i)=P^m$.

Note, in the case where $f_c(P^i)=\left\{\begin{array}{l}P^i \textrm{ for } i <m\\
P^m \textrm{ for }i \geq m\end{array}\right.$, $f_c(0)\subseteq
\bigcap_{i \geq 0} f_c(P^i)=P^m$.  Thus $f_c(0)=(0)$ or $f_c(0)=P^m$
since $f_c(P^n)=P^m$ for $n \geq m$.  Hence, $f_c=f_m$ or $f_c=g_m$
as defined in the statement of the proposition.

Now, suppose that $f_c$ is a semiprime operation which is not
constant on any such interval $m \leq i \leq n$ with $m <n$.
Suppose $f_c(P^i) = P^k$ for $k <i$. Then $P^k=f_c(P^k)\subseteq
f_c(P^j)\subseteq f_c(P^i)=P^k$ for all $k \leq j \leq i$ which
contradicts that fact that $f_c$ is not constant on any interval.
Hence, $f_c(P^i)=P^i$ for all $i\geq 0$.  Since $f_c(0) \subseteq
f_c(P^i)=P^i$ for all $i \geq 0$, then $f_c(0) \subseteq \bigcap_{i
\geq 0} \ P^i=(0)$.  Hence, $f_c$ must be the identity map.

Clearly $f_m \circ f_n= f_{\textrm{min}(m,n)}$ and $g_m \circ g_n=
g_{\textrm{min}(m,n)}$ both imply that the corresponding sets of semiprime
operations in $S_R$, $M_0$ and $M_f$ are submonoids of $M_{\frak{I}}$.
That $M_f$ is a left $M_0$-act can be seen by $f_n \circ g_m =
g_{\textrm{min}(m,n)}$. However, for $m > n$, $g_m \circ f_n(0)=P^m$ and
$(g_m \circ f_n)\circ(g_m \circ f_n)=g_n$ which implies $g_m \circ f_n$ is
not a closure operation. Thus, $M_f$ is not a right $M_0$-act and $S_R =
M_0 \bigcup M_f$ is not a submonoid. \qed

For every $n \geq 0$, $M_n=\{e\} \cup \{f_n\} \cup \{g_n\}$ also
form finite submonoids of $M_{\frak{I}}$ contained in $S_R$,
inter-relating $M_0$ and $M_f$.

\begin{prp}  The only element of $P_R$ when $(R,P)$ is a discrete
valuation ring is the identity.\end{prp}

{\bf Proof:}  Let $(b)=P$.  If $f_c$ is prime, then
$bf_c(P^{i})=f_c(bP^{i})=f_c(P^{i+1})$.  Note if $f_c$ was either
$f_m$ or $g_m$ in the above proof, then $bf_c(P^m)=bP^{m} \subsetneq
P^m=f_c(P^{m+1})=f_c(bP^{m})$.  This contradicts the assumption of
primeness. Thus $P_R=\{e\}$.  \qed

If $R$ is a Dedekind domain which is not necessarily local then for
every maximal ideal $\frak{m}$ in $R$, $R_{\frak{m}}$ is a discrete
valuation ring.  We know the structure $S_{R_{\frak{m}}}$, and can
build the structure of $S_R$ from $S_{R_{\frak{m}}}$.

Given a Dedekind domain $R$ with maximal ideals $P_i$, $i \in
\Lambda$. Consider the monoid given by
$\displaystyle\coprod\limits_{i \in \Lambda} \mathbb{N}_0$, the
coproduct of $\mathbb{N}_0$ (i.e. the set of all functions $\phi:
\Lambda \rightarrow \mathbb{N}_0$ such that $\phi(i)=0$ for all but
finitely many $i \in \Lambda$). Suppose $\phi (i_j)=m_j \neq 0$ for
$i_1,i_2, \ldots i_s$ and $\phi (i)=0$ all other $i$. This $\phi$
corresponds to the ideal $P_{i_1}^{m_1}P_{i_2}^{m_2}\cdots
P_{i_s}^{m_s}$.  The function $\phi \equiv 0$ in $\displaystyle
\coprod\limits_{i \in \Lambda} \mathbb{N}_0$ corresponds to the unit
ideal $R$.

As the non-negative integers play a major role in identifying the
semiprime operations in a discrete valuation, certain subsets of the
semigroup $\mathbb{N}_0^{\Lambda}=\coprod\limits_{i \in \Lambda}
\mathbb{N}_0$ will determine the semiprime operations of a Dedekind
domain with maximal ideals $P_i, i \in \Lambda$. All the nonzero
ideals in a Dedekind domain are finite products of the $P_i$, i.e.
$I=P_{i_1}^{m_1}\cdots P_{i_r}^{m_r}$.

To determine these subsets, first consider the semilocal principal
ideal domain $R$ with two maximal ideals $P$ and $Q$, the ideals of
$R$ are $P^iQ^j$ $i,j \geq 0$ which corresponds to the lattice point
$(i,j)$ in $\mathbb{N}_0^2$. Suppose that for some semiprime
operation $f_c$ defined on $(R,P,Q)$,
$$f_c(P^i)=\left\{ \begin{array}{l} P^i \textrm{ for } i <m\\
P^m \textrm{ for } i \geq m \end{array}\right. \textrm{ and }
f_c(Q^j)=\left\{ \begin{array}{l} Q^j \textrm{ for } j <n\\
Q^n \textrm{ for } j \geq n \end{array}\right..$$

As $f_c$ is semiprime, we know that $$f_c(P^i)f_c(Q^j) \subseteq
f_c(P^iQ^j) \subseteq f_c(P^i) \cap f_c(Q^j)=f_c(P^i)f_c(Q^j)$$ as
$P^iQ^j \subseteq P^i$ and $P^iQ^j \subseteq Q^j.$ Thus

$$f_c(P^iQ^j)= \left\{ \begin{array}{l}P^mQ^n \textrm{ if } i \geq m \textrm{ and } j \geq n\\
                                        P^mQ^j \textrm{ if } i \geq m \textrm{ and } 0 \leq j < n\\
                                        P^iQ^n \textrm{ if } 0 \leq i < m \textrm{ and } j \geq n\\
                                        P^iQ^j \textrm{ if } 0 \leq i < m \textrm{ and } 0 \leq j < n.
                                        \end{array}\right.$$

We define the {\it identity rectangle} $B$ of a semiprime operation
$f_c$ on the lattice $(R,P,Q)$ to be the $(i,j)$ such that
$f_c(P^iQ^j)=P^iQ^j$.

In general, where $\Lambda$ isn't necessarily a two-element set, we
denote the ideal corresponding to $\phi \in \mathbb{N}_0^\Lambda$ by
$I(\phi)$.  Similarly we can define an identity $\Lambda$-box for
$R$ with maximal ideals indexed by $\Lambda$.

\begin{dfn}  The {\it identity $\Lambda$-box} $B_{\Lambda}$ of the semiprime operation $f_c$
over a Dedekind domain $R$ is the set of all $\phi \in
\coprod\limits_{i \in \Lambda} \mathbb{N}_0$ such that
$f_c(I(\phi))=I(\phi).$
\end{dfn}

For simplicity we will denote $\phi_{i}^j$ to be the element of
$\coprod\limits_{i \in \Lambda}\mathbb{N}_0$ such that $\phi(i)=j$ and
$\phi(\lambda)=0$ for all $\lambda \neq i$.  All elements are of the form
$\phi_{i_1}^{j_1}+\phi_{i_2}^{j_2}\cdots +\phi_{i_r}^{j_r}:=\phi_{i_1 i_2
\cdots i_r}^{j_1 j_2 \cdots j_r}$ for distinct $i_k$. Note that the
identity $\Lambda$-box $B_{\Lambda}$ of $f_c$ could be bounded if for
every $i \in \Lambda$ there is a finite $m$ with
$f_c(I(\phi_i^j))=I(\phi_i^m)$ for $j \geq m$.  For each $i_h \in
\Lambda$, define $m_h= \left\{
\begin{array}{l}
m \textrm{ if } f_c(I(\phi_{i_h}^j))=I(\phi_{i_h}^{m}) \textrm{ for }
j \geq m\\
\infty \textrm{ otherwise.}
\end{array}\right.$
In fact, all semiprime operations on the ideals of
$\coprod\limits_{i \in \Lambda}\mathbb{N}_0$ satisfy the equations
$$f_{B_{\Lambda}}(I(\phi_{i_1 i_2 \cdots i_r}^{j_1 j_2 \cdots j_r}))=
\left\{ \begin{array}{l}I(\phi_{i_1 i_2 \cdots i_r}^{j_1 j_2 \cdots j_r})
\textrm{ if } \phi_{i_1,i_2,\ldots,i_r}^{j_1,j_2,\ldots,j_r}
             \in B_{\Lambda}\\
             I(\phi_{i_1 i_2 \cdots i_r}^{k_1 k_2 \cdots k_r}) \textrm{ if } \phi_{i_1,i_2,\ldots,i_r}^{j_1,j_2,\ldots,j_r}
             \notin
             B_{\Lambda} \textrm{ and } k_l=m_l \neq \infty \\
             \qquad \textrm{ for some } l \textrm{ and } k_h=j_h \textrm{ for all } h
             \textrm{ with } k_h \leq m_h.
                                        \end{array}\right.$$

If $B_{\Lambda}$ and $C_{\Lambda}$ are any two identity
$\Lambda$-boxes, clearly, $B_{\Lambda} \cap C_{\Lambda}$ is also a
identity $\Lambda$-box and the action of $f_{B_{\Lambda}} \circ
f_{C_{\Lambda}}$ on nonzero ideals of $R$ is the same as that of
$f_{B_{\Lambda} \cap C_{\Lambda}}.$

Since the semiprime operations of a Dedekind domain correspond to
elements of $\coprod\limits_{i \in \Lambda} \mathbb{N}_0 \cup
\{\infty\}$ under partial ordering, when $B_{\Lambda}$ is bounded
with a finite number $i \in \Lambda$ with $m_i \neq 0$, there are
two types of semiprime operations $f_{B_{\Lambda}}$ and
$g_{B_{\Lambda}}$.  The only difference is that
$f_{B_{\Lambda}}(0)=(0)$ and
$g_{B_{\Lambda}}(0)=P_{i_1}^{m_1}P_{i_2}^{m_2} \cdots
P_{i_r}^{m_r},$ where $\{i_1, \ldots, i_r\}$ is exactly the set of
all $i_j \in \Lambda$ with $m_j < \infty$.

Let us define two subsets of $M_\frak{I}$: \begin{itemize}
\item $M_f=\{ e\} \bigcup\  \{ \ g_{B_\Lambda} \in M_\frak{I}| g_{B_\Lambda}(0)=
P_{i_1}^{m_1}P_{i_2}^{m_2} \cdots P_{i_r}^{m_r} \textrm{ for some
primes } P_{i_j}, j=1,\ldots, r\}$: the set of closure operations
for which the zero ideal is not closed (along with the identity).

\item $M_0=\{ e\} \bigcup\  \{ \ f_{B_\Lambda} \in M_\frak{I}|f_{B_\Lambda}(0)=(0)\}$: the set
of closure operations for which the zero ideal is closed.
\end{itemize}

Suppose now that $B_{\Lambda}$ and $C_{\Lambda}$ are two identity
$\Lambda$-boxes  with both $B_{\Lambda}$ and $C_{\Lambda}$ bounded. Then
$B_{\Lambda} \cap C_{\Lambda}$ is also bounded and is also a identity
$\Lambda$-box and $f_{B_{\Lambda}} \circ f_{C_{\Lambda}}=f_{B_{\Lambda}
\cap C_{\Lambda}}$ and $g_{B_{\Lambda}} \circ
g_{C_{\Lambda}}=g_{B_{\Lambda} \cap C_{\Lambda}}.$ This shows that $M_0$
and $M_f$ are submonoids of $M_\frak{I}$.

Lastly, suppose that $B_{\Lambda}$ and $C_{\Lambda}$ are two identity
$\Lambda$-boxes  with $C_{\Lambda}$ bounded. Then $B_{\Lambda} \cap
C_{\Lambda} \subsetneq C_{\Lambda}$ is also bounded as above and is also a
identity $\Lambda$-box. Note that, $f_{B_{\Lambda}} \circ
g_{C_{\Lambda}}=g_{B_{\Lambda} \cap C_{\Lambda}}$ but $g_{C_{\Lambda}}
\circ f_{B_{\Lambda}}\neq g_{B_{\Lambda} \cap C_{\Lambda}}$ since

$g_{C_{\Lambda}} \circ f_{B_{\Lambda}}(0)=\bigcap\limits_{i \in \Lambda} \
P_i^m $ where $\phi_i^m \in C_{\Lambda} \neq B_{\Lambda} \cap C_{\Lambda}$
which is not a closure operation. This shows that $M_f$ is a left
$M_0$-act, but not a right $M_0$-act.

We have just proved:

\begin{prp}  When $R$ is a Dedekind domain, $S_R$ can
be decomposed into the union of two submonoids $M_0=\{ e\} \bigcup \{
f_{B_\Lambda} \in M_\frak{I}\}$ and $M_f=\{ e\} \bigcup \{ g_{B_\Lambda}
\in M_\frak{I}\}$ where $M_f$ is a left $M_0$-act but not a right
$M_0$-act under composition.\end{prp}

\begin{prp}  The only element of $P_R$ when $R$ is a Dedekind domain is the identity.\end{prp}

{\bf Proof:}  Suppose $(b_i)=P_i$.  If $f_c$ is prime, then
$b_if_c(I)=f_c(b_iI)$ for all $I$. In particular,
$b_if_c(P^j)=f_c(b_i(P_i)^{j})$ for all $j \geq 0$.  Note if $f_c$
was either $f_{B_{\Lambda}}$ or $g_{B_{\Lambda}}$  and $P_i$ is a
prime such that $f_c(P_i^j)=P_i^{m_i}$ for $j \geq m_i$ then
$b_if_c(P_i^{m_i})=P_i^{m_i+1} \subsetneq
P_i^{m_i}=f_c(b_iP_i^{m_i})$. This contradicts the assumption of
primeness. Thus $P_R=\{e\}$. \qed

\section{$S_R$ and $P_R$ when $R=K[[t^2,t^3]]$}

Although $K[[t^2,t^3]]$ is a local ring, the ideal structure in
$K[[t^2,t^3]]$ is not totally ordered as in the case of a discrete
valuation ring.  All ideals in $K[[t^2,t^3]]$ are either generated
by one element $t^n+at^{n+1}$ where $a \in K$ or two elements of the
form $(t^n,t^{n+1})$.  I would like to thank Hwa Young Lee for
pointing out that I was ignoring the ideals $(t^i+at^{i+1})$, with
$a \neq 0$ in a previous version of this paper.  She shared with me
some of the ideas from her developing thesis including some theorems
which she proved which can be summed up in the following
proposition.  The proof here is my own.

\begin{prp}  \label{ideals} Each nonzero nonunit ideal of $R=K[[t^2,t^3]]$ can either be expressed
as a principal ideal in the form $(t^n+at^{n+1})$, $a \in K$, $n
\geq 2$, or as a two generated ideal $(t^n,t^{n+1})$ for $n \geq
2$.\end{prp}

{\bf Proof:}  Suppose $0 \neq f \in R$.  Thus, after multiplying by a
nonzero element of $K$, $f=t^n+a_1t^{n+1}+a_2t^{n+2}+\cdots$ for $n \geq
2$. We will show that $t^{m} \in (f)$ for $m \geq n+2$.  Hence,
$t^n+a_1t^{n+1} \in (f)$. Similarly, $t^{m} \in (t^{n}+a_1t^{n+1})$ for $m
\geq n+2$. Hence, $f \in (t^n+a_1t^{n+1})$.

Let $g \in K[[t]]$. Note that $t^{m-n}g \in K[[t^2,t^3]]$.  Hence, if $g$
is a unit in $K[[t]]$, then $t^{m-n}g^{-1} \in K[[t^2,t^3]]$ also. In
$K[[t]]$, $f=t^{n}(1+a_1t+a_2t^2+\cdots)=t^{n}g$.  Note that
$t^{m-n}g^{-1}f=t^m$.  Similarly $t^m \in (t^{n}+a_1t^{n+1})$.  Since
$f-(t^n+a_1t^{n+1})=a_2t^{n+2}+a_3t^{n+3}+ \cdots \in (f) \bigcap
(t^n+a_1t^{n+1})$, we see that $(t^n+a_1t^{n+1})=(f)$. Hence, all
principal ideals of $K[[t^2,t^3]]$ have the form $(t^n+at^{n+1})$.

Suppose, $I$ is not principal.  As $t^m \in (t^n+at^{n+1})$ for $m \geq
n+2$, then $I$ can be generated by at most 2 elements of the form
$(t^n+at^{n+1},t^m+bt^{m+1})$ where $m=n$ or $m=n+1$.  If $m=n$, then
$t^{n+1} \in I$ which also implies that $t^n \in I$. Hence
$I=(t^n,t^{n+1})$.  If $m=n+1$, then $t^{n+2} \in (t^{n}+at^{n+1})
\subseteq I$ as in the principal case above.  However,
$t^{n+1}=t^{n+1}+bt^{n+2}-bt^{n+2} \in I$ and once again $t^n \in I$.
Hence, $I=(t^n,t^{n+1})$.  \qed

In fact the ideals are woven in the following way:

$$ 
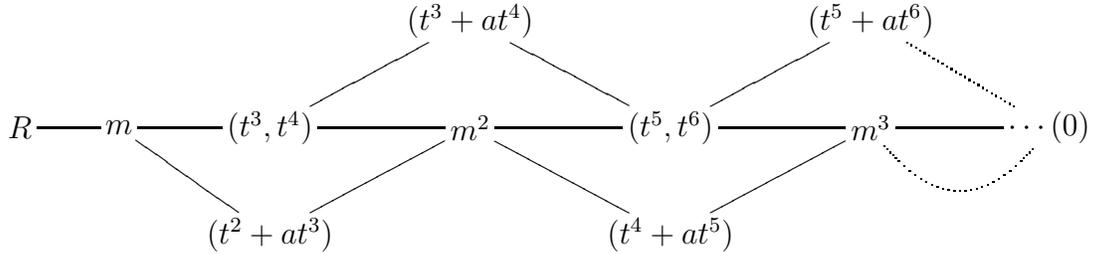
\xymatrix{
                 &       &           & (t^3+at^4) \ar @{-} [dl] \ar @{-} [dr]&          & (t^5+at^6) \ar @{-} [dl]\ar @{.} [dr]&         & \\
    R\ar @{-} [r]& m  \ar @{-} [r]\ar @{- } [dr]& (t^3,t^4)\ar @{-} [r]& m^2\ar @{-} [r] \ar @{-} [dr]&(t^5,t^6)\ar @{-} [r]&m^3 \ar @{-} [r]&\cdots  (0)\ar @/^2pc/ @{.} [l]\\
                 &           & (t^2+at^3) \ar @{-} [ur]     &       &  (t^4+at^5) \ar @{-} [ur]  &       &       } $$
where each line segment in the above diagram indicates $\supseteq$.

In the case of a discrete valuation ring $(R,P)$, integral closure
is the identity map on ideals of $R$.  For $K[[t^2,t^3]]$, the
integral closure of ideals of the form $(t^i+at^{i+1})$ is
$\overline{(t^i+at^{i+1})}=(t^i,t^{i+1})$ whereas the ideals of the
form $(t^i,t^{i+1})$ are all integrally closed.  Looking at the
above diagram, we see that the chain of ideals in the center are all
integrally closed. However, the principal ideals are not. Clearly
there are now many more closure operations for $K[[t^2,t^3]]$.  In
fact, the semiprime operations which are not bounded abound.  To
shorten the expressions appearing in the proofs we will denote the
principal ideals $P_{i,a}:=(t^i+at^{i+1})$ and
$M_{i}:=(t^i,t^{i+1})$.

\begin{prp}  In $K[[t^2,t^3]]$, for all $i \geq 2$ and all $a \in K$, the map \newline
$f_i^{int}(I)
:=\left\{ \begin{array}{l} M_i \textrm{ if } I=P_{i,a} \\
              I \textrm{ if } I \neq P_{i,a} \end{array}\right.$ is a closure
operation which is not semiprime.
\end{prp}

{\bf Proof:}  Clearly $f_i^{int}(I) \supseteq I$ for all $I$ and if $I
\subseteq J$, $f_i^{int}(I) \subseteq f_i^{int}(J)$.  As $f_i^{int}(I)=I$
whenever $I \neq P_{i,a}$, and $$f_i^{int}\circ
f_i^{int}(P_{i,a})=f_i^{int}(P_{i,a})=M_i=f_i^{int}(M_i)$$ then
$f_i^{int}$ is a closure operation.

As $M_jP_{k,a}=M_{j+k}$, the only ideals which are proper factors of
$P_{m,a}$ are of the form $P_{j,b}$, $j\leq m-2$ and $b \in K$. If $j+k=m$
with $j,k \geq 2$, then
$$f_i^{int}(P_{j,a})f_i^{int}(P_{k,b})=
\left\{ \begin{array}{l}P_{m,a+b} \textrm{ if } j \neq i \textrm{ and } k \neq i \\
              M_m \textrm{ if } j=i \textrm{ or } k=i. \end{array}\right.$$
If $m\geq i+j$, $j \geq 2$, $f_i^{int}(P_{m,a+b})=P_{m,a+b}$ and
$M_m \nsubseteq P_{m,a+b}$. Thus $f_i^{int}$ is not a semiprime
operation. \qed

We observe in the proof, that if we want such a closure operation
which maps $P_{i,a}$ to $M_i$ to be semiprime we also need $P_{m,a}$
to map to $M_m$ for $m \geq i+2$. Hence, we have the following:

\begin{crl}  Let $S \neq \emptyset$ and $T$, possibly empty, be subsets of the field $K$.
Over $K[[t^2,t^3]]$ for all $i \geq 2$, the maps $$f_{i,S,T}^{int}(I)=
\left\{ \begin{array}{l} I \textrm{ if } I \supseteq M_{i+1}, I=P_{i,a}, a \notin S \textrm{ or } I = P_{i+1,b}, b \notin T\\
              \overline{I} \textrm{ if } I \subseteq M_{i+2}, I \supseteq P_{i,a}, a \in S \textrm{ or } I \subseteq P_{i+1,b}, b \in
              T\end{array}\right.$$
are semiprime operations.
\end{crl}

{\bf Proof:}  Clearly $f_{i,S,T}^{int}$  are also closure operations
and from the proof of above, they are semiprime. \qed

\begin{lemma} \label{plus2} If $f_c$ is a semiprime operation on $K[[t^2,t^3]]$
and $M_j=f_c(M_{j+2})$ for some $j$, then $f_c$ is a bounded
semiprime operation.\end{lemma}

{\bf Proof:}  As $M_j \supseteq M_{j+1} \supseteq M_{j+2},$ then
$$M_j=f_c(M_{j+2})\subseteq f_c(M_{j+1})\subseteq f_c(M_j)=f_c(f_c(M_{j+2}))=M_j.$$

We will use induction to show that $f_c(M_{j+n})=M_j$ for $n \geq 0$.
Assume that $f_c(M_{j+k})=M_j$ for $2 \leq k \leq n$.  Since
$M_{j+n+1}=P_{2,0}M_{j+n-1}$,
$$f_c(M_{j+n+1}) \supseteq
f_c(P_{2,0})f_c(M_{j+n-1}) \supseteq P_{2,0}f_c(M_{j+k-2}) = M_{j+2}
\supseteq M_{j+n+1}.$$  Applying $f_c$ to the chain,
$f_c(M_{j+n+1})\supseteq f_c( M_{j+2})=M_j \supseteq f_c(M_{j+n+1}).$  As
the right hand and left hand sides of the chains are equal, we obtain
$f_c(M_{j+n+1})=M_j$.

For any $a \in K$ and $k \geq 0$, we have $M_{j+k} \supseteq P_{j+k,a}
\supseteq M_{j+k+2}.$ Applying $f_c$ to the chain and using the fact that
$f_c(M_{j+k})=M_j$ for $k \geq 0$ we obtain $f_c(P_{j+k,a})=M_j.$

Since the above arguments show if $0 \neq I \subseteq M_j$,
$f_c(I)=M_j$, by the definition of bounded, we see that $f_c$ is a
bounded semiprime operation. \qed

\begin{lemma} \label{plus2'} If $f_c$ is a semiprime operation on $K[[t^2,t^3]]$
and $f_c(M_j)=f_c(M_{j+2})$ for some $j$, then $f_c$ is a bounded
semiprime operation.\end{lemma}

{\bf Proof:}  We can break the proof down into the following two
cases:

\begin{enumerate}
\item $f_c(M_j)=M_k$, $k \leq j$ or
\item $f_c(M_j)=P_{k,a}$, some $a \in K$ and $k \leq j-2$
\end{enumerate}

In case (1), $M_k=f_c(M_k) \supseteq f_c(M_{k+1}) \supseteq f_c(M_{k+2})
\supseteq f_c(M_{j+2})=M_k.$ By Lemma \ref{plus2}, $f_c$ is bounded.

In case (2), we need to show that for any nonzero ideal $I \subseteq
P_{k,a}$, $f_c(I)=P_{k,a}$.  Clearly, if $M_{j+2} \subseteq I
\subseteq P_{k,a}$, then $P_{k,a}= f_c(M_j)=f_c(M_{j+2})=f_c(I).$ We
will see by induction that $f_c(M_{j+n})=P_{k,a}$ for $n \geq 2$.
Assume that $f_c(M_{j+i})=P_{k,a}$ for $2 \leq i \leq n$.  Since
$M_{j+n+1}=P_{2,0}M_{j+n-1}$, we have
$$f_c(M_{j+n+1})=f_c(P_{2,0}M_{j+n-1}) \supseteq
f_c(P_{2,0})f_c(M_{j+n-1}) \supseteq P_{k+2,a} \supseteq M_{j+n+1}.$$ Note
that $M_{j+2} \subseteq P_{k+2,a} \subseteq P_{k,a}$.  Hence,
$f_c(P_{k+2,a})=P_{k,a}$ which implies after applying $f_c$ to the above
chain of containments that $f_c(M_{j+n+1})=P_{k,a}$.  Hence,
$f_c(M_{j+n})=P_{k,a}$ for $n \geq 0$.

Since $M_{k+n} \supseteq P_{k+n,b} \supseteq M_{k+n+2}$, applying
$f_c$ to this chain of containments and noting that
$f_c(M_{k+n})=P_{k,a}$ for all $n \geq 2$, we conclude that
$f_c(P_{k+n,b})=P_{k,a}$. Now we have seen that for all nonzero $I
\subseteq P_{k,a}$,  $f_c(I)=P_{k,a}$. Hence, $f_c$ is bounded. \qed

\begin{lemma} \label{plus1} If $f_c$ is a semiprime operation on $K[[t^2,t^3]]$
and $f_c(M_j)=f_c(M_{j+1})$ for some $j$, then $f_c$ is a bounded
semiprime operation.\end{lemma}

{\bf Proof:}  Note that for $j \geq 2$, if $R=f_c(M_j)$ then
$f_c(M_{2j})=f_c(M_{j}^2)\supseteq f_c(M_{j})^2=R.$  Since $M_j
\supseteq M_{j+2}\supseteq M_{2j}$, then
$R=f_c(M_j)=f_c(M_{j+2})=f_c(M_{2j})$. By Lemma \ref{plus2'} we can
conclude that $f_c$ is bounded.

Also for $j \geq 3$ if $I=f_c(M_j) \supseteq M_{j-1} \supseteq M_j$
then $I=f_c(M_{j-1})=f_c(M_{j+1}).$ By Lemma \ref{plus2'} we can
conclude that $f_c$ is bounded.  That leaves us with the cases:

\begin{enumerate}
\item $f_c(M_j)=M_j$ for $j \geq 2$ or
\item $f_c(M_j)=P_{j-2,a}$ for some $a \in K$ and $j \geq 4$.
\end{enumerate}

In case (1), consider $f_c(M_{2j+2}) \supseteq
f_c(M_{j+1}^2)\supseteq f_c(M_{j+1})^2=M_{2j} \supsetneq M_{2j+2}.$
Applying $f_c$, we now see that $f_c(M_{2j})=f_c(M_{2j+2}).$  Again,
Lemma \ref{plus2'} yields that $f_c$ is bounded.

In case (2), $f_c(M_{j-1}) \supseteq f_c(M_j)=P_{j-2,a} \textrm{ and }
f_c(M_{j-1}) \supseteq M_{j-1}.$  Thus
$$f_c(M_{j-1})\supseteq P_{j-2,a}+M_{j-1}=M_{j-2} \supseteq
M_{j-1} \textrm{ implies } f_c(M_{j-1})=f_c(M_{j-2}).$$ Now we are in the
same set up as our Lemma but two steps up.  If $f_c(M_{j-2})=M_{j-2}$ we
are done by case (1) above. Otherwise,
$f_c(M_{j-2})=f_c(M_{j-1})=P_{j-4,b}$, some $b \in K$.  Now

Now $P_{2j-8,2b}=f_c(M_{j-1})^2 \subseteq f_c(M_{2j-2}) \subseteq
f_c(M_{2j-4}) \subseteq f_c(P_{2j-8,2b})$.  If we apply $f_c$ to
this chain of containments we see that
$f_c(M_{2j-2})=f_c(M_{2j-4}).$   Again $f_c$ is bounded by Lemma
\ref{plus2'}. \qed

\begin{lemma} \label{plus1'} If $f_c$ is a semiprime operation on $K[[t^2,t^3]]$
and $f_c(M_j)=f_c(P_{j-2,b})$ for some $j \geq 4$ and $b \in K$,
then $f_c$ is a bounded semiprime operation.

\end{lemma}

{\bf Proof:}  As in the proof of Lemma \ref{plus1}
$f_c(M_j)=f_c(P_{j-2,b})$,  implies that
$f_c(M_{j-1})=f_c(M_{j-2})$.  We now conclude by Lemma \ref{plus1}
that $f_c$ is also bounded. \qed

The following theorem describes the unbounded semiprime operations
over $K[[t^2,t^3]]$.

\begin{thrm}  \label{usp} Let $S$ be a nonempty subset of $K$, $T$ any subset.
If $f_c$ is an unbounded semiprime operation over $K[[t^2,t^3]]$,
then $f_c$ is either the identity or
$$f_{i,S,T}^{int}(I)
\left\{ \begin{array}{l} I \textrm{ if } I \supseteq P_{i,a}, a \notin S \textrm{ or } I \supseteq P_{i+1,b}, b \notin T\\
              \overline{I} \textrm{ if } I \subseteq M_{i+2}, I=P_{i,a}, a \in S \textrm{ or } I = P_{i+1,b}, b \in T
              \end{array}\right..$$
              \end{thrm}

{\bf Proof:}  Suppose $f_c$ is an unbounded semiprime operation over
$K[[t^2,t^3]]$ which is not the identity. Then $f_c(I) \neq I$ for
some nonzero ideal $I$.

If $I=M_j$ for some $j \geq 2$, then by Lemmas \ref{plus2'}, \ref{plus1}
and \ref{plus1'}, $f_c$ would be bounded, contradicting the unbounded
assumption.  Thus $I$ must be a principal ideal.

If $f_c(P_{k,a})=f_c(P_{k+2,b})$ for some $k$, then
$f_c(P_{k,a})=f_c(M_{k+2})=f_c(P_{k+2,b})$ and Lemma \ref{plus1'} implies
that $f_c$ is bounded, contradicting the unboundedness assumption.

If $f_c(P_{k,a})=f_c(M_{k-1})$ for some $k$, then
$f_c(P_{k,a})=f_c(M_{k})=f_c(M_{k-1})$ which is bounded by Lemma
\ref{plus1}.

Thus $f_c(P_{k,a})=M_k$.  Let $W=\{k|f_c(P_{k,a})=M_k \textrm{ for some }
a \in K  \textrm{ and some }k \geq 2 \}$. Since $W$ is nonempty subset of
the positive integers there is a smallest $j \geq 2$ in $W$.  Let $S=\{a
\in K |f_c(P_{j,a})=M_j\}$. Since $P_{n,b}=P_{j,a}P_{n-j,b-a}$ for all $b
\in K$, $a \in S$ and all $n \geq j+2$, then $f_c(P_{n,b}) \supseteq
f_c(P_{j,a})f_c(P_{n-j,b-a})\supseteq M_jP_{n-j,b-a}=M_n \supseteq
P_{n,b}.$  Applying $f_c$ to the chain of containments, we see that
$M_n=f_c(P_{n,b})$ for all $b \in K$ and $n=j$ or $n \geq j+2$.

Note for all $b \in K$, $f_c(P_{j+1,b}) \subseteq f_c(M_{j+1})=M_{j+1},$
thus $f_c(P_{j+1,b})=P_{j+1,b}$ or $f_c(P_{j+1,b})=M_{j+1}$. If $T= \{b
\in K|f_c(P_{j+1,b})=M_{j+1}\}$ then $f_c=f_{j,S,T}^{int}$ as defined in
the statement of the theorem.\qed

The bounded semiprime operations are given by the following theorem:

\begin{thrm} \label{bsp} The only bounded semiprime operations on $K[[t^2,t^3]]$
are of the forms $$f_{m,a}^f(I)=\left\{ \begin{array}{l}I \textrm{ for } I \supseteq P_{m,a}\\
                         M_{m-1} \textrm{ for } I=P_{m-1,b}, \forall b \in K\\
                         M_{m} \textrm{ for } I=P_{m,b}, b\neq a, I=P_{m+1,d}, \forall d \in K \textrm{ or } I=M_{m+1} \\
                         P_{m,a} \textrm{ for nonzero } I \subseteq
                         P_{m,a}\\
                         (0) \textrm{ if } I=(0)\\
                         \end{array}\right.$$

$$g_{m,a}^f(I)=\left\{ \begin{array}{l}I \textrm{ for } I \supseteq P_{m,a}\\
                         M_{m-1} \textrm{ for } I=P_{m-1,b}, \forall b \in K\\
                         M_{m} \textrm{ for } I=P_{m,b}, b\neq a, I=P_{m+1,d}, \forall d \in K \textrm{ or } I=M_{m+1} \\
                         P_{m,a} \textrm{ for } I \subseteq
                         P_{m,a}\\
                         \end{array}\right.$$
for $m \geq 2$ and $a \in K$,

$$f_{n,S,T,m}^f(I)=\left\{ \begin{array}{l}I \textrm{ for } I \supseteq P_{n,a}, a \notin S \textrm{ or } I \supseteq P_{n+1,b}, b \notin T\\
                         \overline{I} \textrm{ for }  P_{m-1,d} \subseteq I \subseteq J, J=P_{n,a}, a \in S \textrm{ or }
                         J=P_{n+1,b}, b \in T \\
                         M_m \textrm{ for nonzero } {I} \subseteq M_m\\
                         (0) \textrm{ if } I=(0)\\
                         \end{array}\right.$$

$$g_{n,S,T,m}^f(I)=\left\{ \begin{array}{l}I \textrm{ for } I \supseteq P_{n,a}, a \notin S \textrm{ or } I \supseteq P_{n+1,b}, b \notin T\\
                         \overline{I} \textrm{ for }  P_{m-1,d} \subseteq I \subseteq J, J=P_{n,a}, a \in S \textrm{ or }
                         J=P_{n+1,b}, b \in T \\
                         M_m \textrm{ for } {I} \subseteq M_m\\
                         \end{array}\right.$$
for $m-1 \geq n \geq 2$, $S \neq \emptyset$ and if $m=n+1$, $T=K$,

$$f_{n,S,T,m\prime}^f(I)=\left\{ \begin{array}{l}I \textrm{ for } I \supseteq P_{n,a}, a \notin S \textrm{ or } I \supseteq P_{n+1,b}, b \notin T\\
                         \overline{I} \textrm{ for }  M_{m-2} \subseteq I \subseteq J, J=P_{n,a}, a \in S \textrm{ or }
                         J=P_{n+1,b}, a \in S, b \in T \\
                         M_{m-2} \textrm{ for } I=P_{m-1,d},M_{m-1}\\
                         M_m \textrm{ for nonzero } {I} \subseteq M_m\\
                         (0) \textrm{ if } I=(0)\\
                         \end{array}\right.$$

$$g_{n,S,T,m\prime}^f(I)=\left\{ \begin{array}{l}I \textrm{ for } I \supseteq P_{n,a}, a \notin S \textrm{ or } I \supseteq P_{n+1,b}, b \notin T\\
                         \overline{I} \textrm{ for }  M_{m-2} \subseteq I \subseteq J, J=P_{n,a}, a \in S \textrm{ or }
                         J=P_{n+1,b}, a \in S, b \in T \\
                         M_{m-2} \textrm{ for } I=P_{m-1,d},M_{m-1}\\
                         M_m \textrm{ for } {I} \subseteq M_m\\
                         \end{array}\right.$$
for $m-2 \geq n \geq 2$, $S \neq \emptyset$.

\end{thrm}

{\bf Proof:}  If $f_c$ is a bounded semiprime operation, then for
small nonzero $I$, either
\begin{enumerate}
\item $f_c(I) =P_{m,a}$ for $I \subseteq P_{m,a}$ some $a \in K$ or
\item $f_c(I)=M_m$ for $I \subseteq M_m$
\end{enumerate} for some $m \geq 2$.

Case (1):  If $f_c(I)=P_{m,a}$ for $I \subseteq P_{m,a}$, then for
$f_c$ to be semiprime, we see as in the proof of Proposition
\ref{ideals} that $f_c(P_{i,b})=P_{i,b}$ for $2 \leq i \leq m-2$ and
all $b \in K$, since the only factors of $P_{m,a}$ are $P_{i,b}$ for
$2 \leq i \leq m-2$ any $b \in K$.  Note that $f_c(M_i) \subseteq
f_c(P_{i-2,b})=P_{i-2,b}$ for all $4 \leq i \leq m-2$ and $b \in K$.
Hence, $f_c(M_i) \subseteq \bigcap\limits_{b \in K} P_{i-2,b}=M_i$
for $4 \leq i \leq m-2$.  $M_2 \supseteq M_3$ contain only the unit
principal ideal $R$.  Let $i=2,3$, then $M_jf_c(M_i) \subseteq
f_c(M_i)f_c(M_j) \subseteq f_c(M_{i+j})=M_{i+j}$, for $2 \leq j \leq
m-2-i$.  This set of containments implies that  $f_c(M_i) \subseteq
M_{i+j}:M_j=M_i$.  Hence $f_c(M_i)=M_i$ for $i=2,3$.

Since $M_{m+3} \subseteq P_{m,a}$, we see that
$f_c(M_{m+3})=P_{m,a}$. Applying $f_c$ to the following chain of
containments: $M_{m+3} \subseteq P_{m+1,b} \subseteq M_{m+1}$, we
see that $P_{m,a} \subseteq f_c(P_{m+1,b})\subseteq f_c(M_{m+1}).$
However, $P_{m+1,b}\subseteq M_{m+1}$ are both incomparable with
$P_{m,a}$. Thus $M_m \subseteq f_c(P_{m+1,b}) \subseteq
f_c(M_{m+1}).$  Since $P_{2,d}P_{m+1,b}=P_{m+3,b+d}$ for all $d \in
K$, we see that
$$P_{2,d}f_c(P_{m+1,b}) =f_c(P_{2,d})f_c(P_{m+1,b})\subseteq
f_c(P_{m+3,b+d})=P_{m,a}.$$ Thus $f_c(P_{m+1,b}) \subseteq P_{m-2,a-d} $
for all $d \in K$. Since $\bigcap\limits_{d \in K} P_{m-2,a-d}=M_m$, we
see that $f_c(P_{m+1,b})=M_{m}.$  Now since $P_{m+1,b} \subseteq M_{m+1}
\subseteq M_m$ then we easily see that $f_c(M_{m+1})=M_m$ also.

As $M_{m+1} \subseteq P_{m-1,d}$ for all $d \in K$ and
$f_c(M_{m+1})=M_{m}$, we can see that $M_{m-1}=M_m+P_{m-1,d} \subseteq
f_c(P_{m-1,d}) \subseteq f_c(M_{m-1})$.  Applying $f_c$ to the chain of
containments, we observe that $f_c(P_{m-1,d})=f_c(M_{m-1}).$   Noting that
$f_c(M_{m-1}) \subseteq f_c(P_{m-3,d})=P_{m-3,d}$ for all $d \in K$, we
can conclude that $f_c(M_{m-1}) \subseteq \bigcap\limits_{d \in K}
P_{m-3,d}=M_{m-1}$.  Putting this fact together with the equality above we
see that $f_c(P_{m-1,d})=f_c(M_{m-1})=M_{m-1}.$

Putting all the above arguments together we see that $f_c$ must be
$f_{m,a}^f$ or $g_{m,a}^f$ depending on whether or not $f_c(0)$ is $(0)$
or $P_{m,a}$.  The following diagram represents $f_c$.  The arrows
represent the $f_c$-closure of the indicated ideals.

$$ 
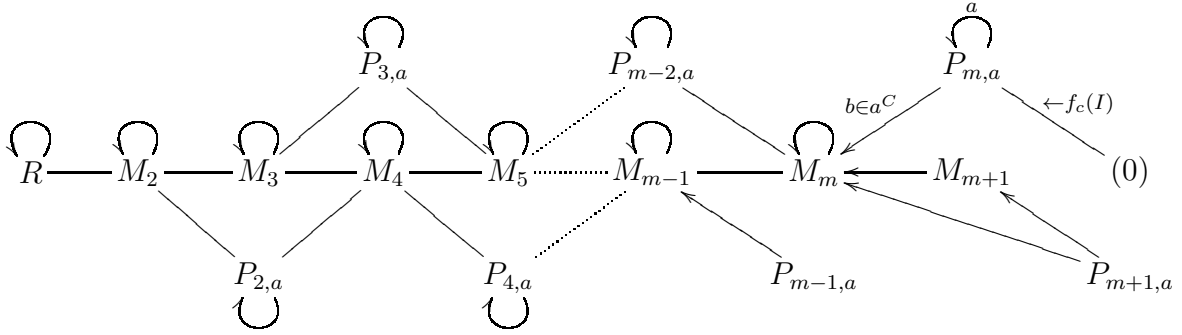
\xymatrix{
                 &       &           & P_{3,a} \ar @(ru,lu)[]\ar @{-} [dl] \ar @{-} [dr]&          & P_{m-2,a}\ar @(ru,lu)[] \ar @{.} [dl]\ar @{-} [dr]&         & P_{m,a}\ar @(ru,lu)[]_{a} \ar @{->}_{b \in a^C} [dl] \ar @{-}^{\leftarrow f_c(I)} [dr]& \\
    R\ar @(ru,lu)[] \ar @{-} [r]& M_2\ar @(ru,lu)[]  \ar @{-} [r]\ar @{- } [dr]& M_3\ar @(ru,lu)[] \ar @{-} [r]& M_4\ar @(ru,lu)[] \ar @{-} [r] \ar @{-} [dr]&M_{5}\ar @(ru,lu)[] \ar @{.} [r]&M_{m-1}\ar @(ru,lu)[] \ar @{-} [r]&M_{m}\ar @(ru,lu)[] &M_{m+1} \ar @{->} [l]&(0)\\
                 &           & P_{2,a}\ar @(rd,ld)[] \ar @{-} [ur]     &       &  P_{4,a}\ar @(rd,ld)[] \ar @{.} [ur]  &       &P_{m-1,a} \ar @{->} [ul]     &  &P_{m+1,a} \ar @{->}[ul]\ar @{->}[ull] } $$

Case (2):  Suppose $f_c(I)=M_m$ for all $I \subseteq M_{m}$ and $m \geq
2$. The closure of ideals in the following diagram still needs to be
determined.

$$ 
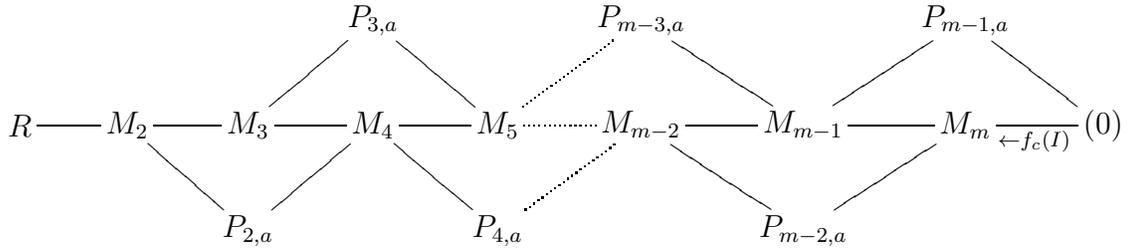
\xymatrix{
                 &       &           & P_{3,a} \ar @{-} [dl] \ar @{-} [dr]&          & P_{m-3,a} \ar @{.} [dl]\ar @{-} [dr]&         & P_{m-1,a} \ar @{-} [dl] \ar @{-} [dr]& \\
    R\ar @{-} [r]& M_2  \ar @{-} [r]\ar @{- } [dr]& M_3\ar @{-} [r]& M_4\ar @{-} [r] \ar @{-} [dr]&M_{5}\ar @{.} [r]&M_{m-2} \ar @{-} [r]&M_{m-1} &M_m \ar @{-} [l]&(0)\ar @{-}^{\leftarrow f_c(I)} [l]\\
                 &           & P_{2,a} \ar @{-} [ur]     &       &  P_{4,a} \ar @{.} [ur]  &       &P_{m-2,a} \ar @{-} [ul]  \ar @{-} [ur]   &  & } $$
Only the $P_{m-1,b}$ are not comparable to $M_m$. Since $M_{m+1} \subseteq
P_{m-1,b} \textrm{ and } f_c(M_{m+1})=M_{m}$ then
$f_c(P_{m-1,b})=f_c(M_{m-1}).$  We will get back to this later; however,
in the next diagram we will indicate this with an arrow from the
$P_{m-1,a}$'s to $M_{m-1}$ and omitting the line from $(0)$ to the
$P_{m-1,a}$'s.

First, we will see that $f_c(M_2)=M_2$. Suppose $f_c(M_2)=R$, then
$f_c(M_{m-2})=f_c(M_2)f_c(M_{m-2}) \subseteq f_c(M_m)=M_m \subseteq
M_{m-2}$.  Applying $f_c$ to this chain of containments, we see that
$M_m=f_c(M_m)=f_c(M_{m-2})$ which is a contradiction since $M_{m-2}
\nsubseteq M_m$.

We now show that $f_c(M_{n})=M_{n}$ for $2 < n \leq m-2$. Suppose $M_n
\subsetneq f_c(M_n)=I$ where $I \supseteq M_{n-1}$ or $I \supseteq
P_{n-2,d}$ for some $d$.  Once again, we decompose $M_m=M_{n}M_{m-n}$.
Noting that $P_{n-2,d}M_{m-n}=M_{m-2}$ and $M_{n-1}M_{m-n}=M_{m-1}$ and
$M_{m-1} \subseteq M_{m-2}$, we see that $M_{m-1} \subseteq
f_c(M_n)f_c(M_{m-n}) \subseteq f_c(M_n)f_c(M_{m-n}) \subseteq f_c(M_m)=M_m
\subseteq M_{m-1}$. As above, this implies that
$M_m=f_c(M_m)=f_c(M_{m-1})$ which gives a contradiction.  These arguments
imply that all ideals along the central line in the above figure excluding
possibly $M_{m-1}$ are $f_c$-closed which I will indicate by a loop in the
diagram.

$$ 
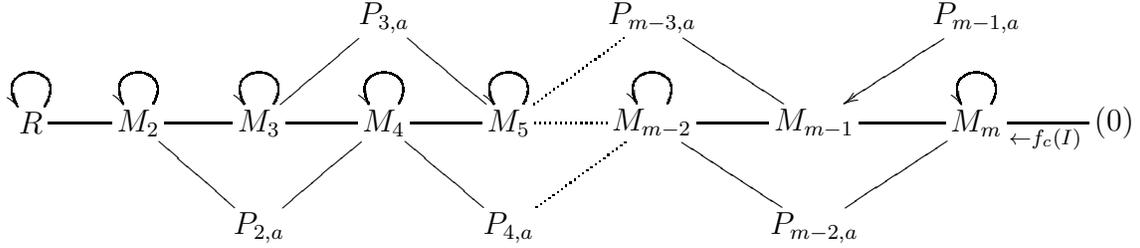
\xymatrix{
                 &       &           & P_{3,a} \ar @{-} [dl] \ar @{-} [dr]&          & P_{m-3,a} \ar @{.} [dl]\ar @{-} [dr]&         & P_{m-1,a} \ar @{->} [dl] & \\
    R\ar @(ru,lu)[] \ar @{-} [r]& M_2\ar @(ru,lu)[]  \ar @{-} [r]\ar @{- } [dr]& M_3\ar @(ru,lu)[]\ar @{-} [r]& M_4\ar @(ru,lu)[]\ar @{-} [r] \ar @{-} [dr]&M_{5}\ar @(ru,lu)[]\ar @{.} [r]&M_{m-2}\ar @(ru,lu)[] \ar @{-} [r]&M_{m-1} &M_m\ar @(ru,lu)[] \ar @{-} [l]&(0)\ar @{-}^{\leftarrow f_c(I)} [l]\\
                 &           & P_{2,a} \ar @{-} [ur]     &       &  P_{4,a} \ar @{.} [ur]  &       &P_{m-2,a} \ar @{-} [ul]  \ar @{-} [ur]   &  & } $$

Now, we will determine $f_c(P_{k,b})$ for $2 \leq k \leq m-2$.  Since
$f_c(P_{k,b}) \subseteq f_c(M_k)=M_k$ for all $b \in K$ and $2 \leq k \leq
m-2$, we see that $f_c(P_{k,b})$ may equal $P_{k,b}$ or $M_k$.  Suppose
that $f_c(P_{k,b})=M_k$ for some $b \in K$ with $2 \leq k \leq m-2$.
Assume $n$ is the smallest $2 \leq n \leq m-2$ satisfying this property
for some $b\in K$ and define $S=\{b \in K\mid f_c(P_{n,b})=M_n\}$.
Observing that $P_{k,d}=P_{j,a}P_{k-j,b}, a+b=d$ for $2 \leq j \leq k-n$
and $f_c(P_{k,d})\supseteq f_c(P_{j,a})f_c(P_{k-j,b}) =M_k$, we conclude
that $f_c(P_{k,d})=M_k$ for all $d \in K$ and $n+2 \leq k \leq m-2$.  For
each of these $P_{k,d}$'s, we indicate that the closure is $M_k$ in the
following diagram by indicating an arrow between $P_{k,d}$ and $M_k$ and
omitting the line between the $P_{k,d}$'s and $M_{k+2}$.  I have left off
the ideals containing $M_{n-1}$ since all of these ideals are now known to
be $f_c$-closed.

$$ 
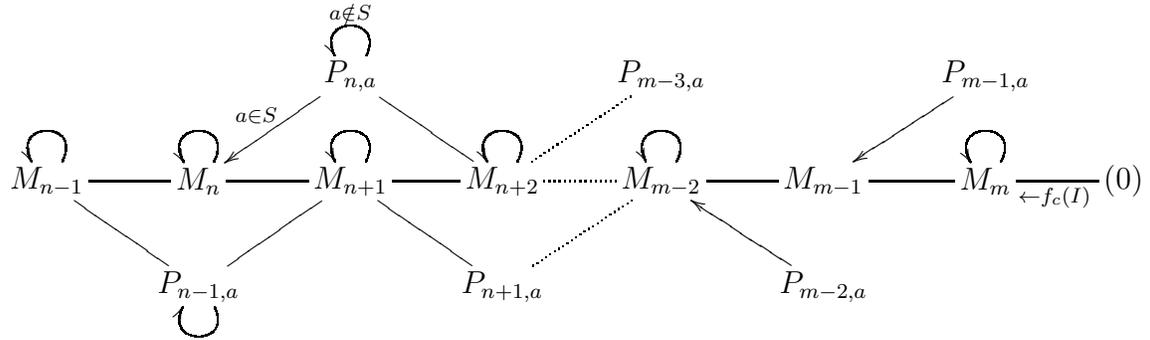
\xymatrix{
                        &          & P_{n,a}\ar @(ru,lu)[]_{a \notin S} \ar @{->}_{a \in S} [dl] \ar @{-} [dr]&          & P_{m-3,a} \ar @{.} [dl]&         & P_{m-1,a} \ar @{->} [dl] & \\
     M_{n-1} \ar @(ru,lu)[] \ar @{-} [r] \ar @{- } [dr]& M_{n}\ar @(ru,lu)[]\ar @{-} [r]& M_{n+1}\ar @(ru,lu)[]\ar @{-} [r] \ar @{-} [dr]&M_{n+2}\ar @(ru,lu)[]\ar @{.} [r]&M_{m-2}\ar @(ru,lu)[] \ar @{-} [r]&M_{m-1} &M_m\ar @(ru,lu)[] \ar @{-} [l]&(0)\ar @{-}^{\leftarrow f_c(I)} [l]\\
               & P_{n-1,a} \ar @(rd,ld)[] \ar @{-} [ur]     &      &  P_{n+1,a} \ar @{.} [ur]  &       &P_{m-2,a} \ar @{->} [ul]     &  & } $$

At this point, there are two ambiguities.  What is $f_c(M_{m-1})$ and what
is $f_c(P_{n+1,b})$ for $b \in K$?  Since $f_c(M_{m-1}) \subseteq
f_c(M_{m-2})=M_{m-2}$ and $f_c(M_{m-1}) \subseteq f_c(P_{m-3,b})$ for all
$b \in K$.  If $f_c(P_{m-3,d})=P_{m-3,d}$ for some $d \in K$, then
$f_c(M_{m-1})\subseteq M_{m-2} \bigcap P_{m-3,d}=M_{m-1}$.  Hence,
$f_c(M_{m-1})=f_c(P_{m-1,b})=M_{m-1}$.  Otherwise,
$f_c(P_{m-1,b})=f_c(M_{m-1})$ could be $M_{m-1}$ or $M_{m-2}$.

In the case that $f_c(P_{m-1,b})=f_c(M_{m-1})=M_{m-1}$, let $S=\{b \in
K|f_c(P_{n,b})=M_n\}$ as above and $T=\{b \in K|f_c(P_{n+1,b})=M_{n+1}\}$,
then $f_c=f_{n,S,T,m}^f$ or $f_c=g_{n,S,T,m}^f$ depending on where $f_c$
maps $(0)$.  From the previous diagram, I have added the loop at $M_{m-1}$
to indicate that $M_{m-1}$ is $f_c$-closed.

$$ 
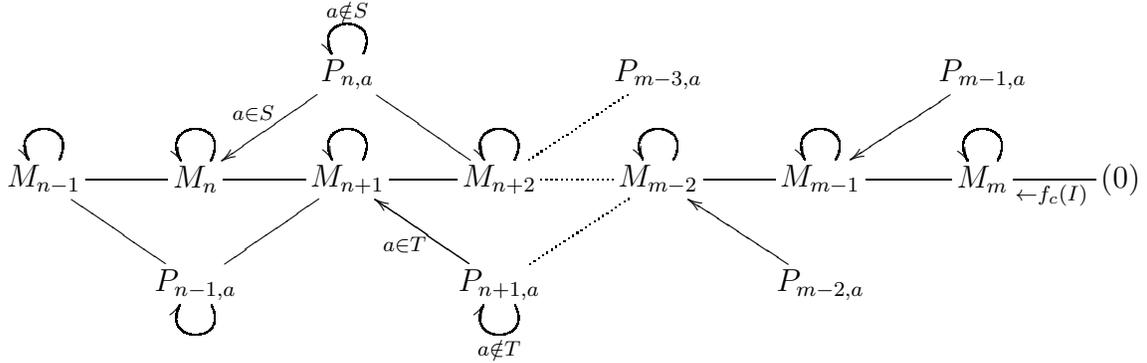
\xymatrix{
                        &          & P_{n,a}\ar @(ru,lu)[]_{a\notin S} \ar @{->}_{a \in S} [dl] \ar @{-} [dr]&          & P_{m-3,a} \ar @{.} [dl]&         & P_{m-1,a} \ar @{->} [dl] & \\
     M_{n-1} \ar @(ru,lu)[] \ar @{-} [r] \ar @{- } [dr] & M_{n}\ar @(ru,lu)[]\ar @{-} [r]& M_{n+1}\ar @(ru,lu)[]\ar @{-} [r] \ar @{-} [dr]&M_{n+2}\ar @(ru,lu)[]\ar @{.} [r]&M_{m-2}\ar @(ru,lu)[] \ar @{-} [r]&M_{m-1}\ar @(ru,lu)[] &M_m\ar @(ru,lu)[] \ar @{-} [l]&(0)\ar @{-}^{\leftarrow f_c(I)} [l]\\
               & P_{n-1,a}\ar @(rd,ld)[] \ar @{-} [ur]     &      &  P_{n+1,a}\ar @(rd,ld)[]^{a \notin T} \ar @{->}^{a \in T} [ul] \ar @{.} [ur]  &       &P_{m-2,a} \ar @{->} [ul]     &  & } $$

In the case that $f_c(P_{m-1,b})=f_c(M_{m-1})=M_{m-2}$, let $S=\{b \in
K|f_c(P_{n,b})=M_n\}$ as above and $T=\{b \in K|f_c(P_{n+1,b})=M_{n+1}\}$,
then $f_c=f_{n,S,T,m\prime}^f$ or $f_c=g_{n,S,T,m\prime}^f$ depending on
where $f_c$ maps $(0)$.  Unlike the previous diagram, there is not a loop
at $M_{m-2}$, since it is not $f_c$-closed, but an arrow from both
$M_{m-1}$ and the $P_{m-1,a}$'s to indicate their $f_c$-closure.

$$ 
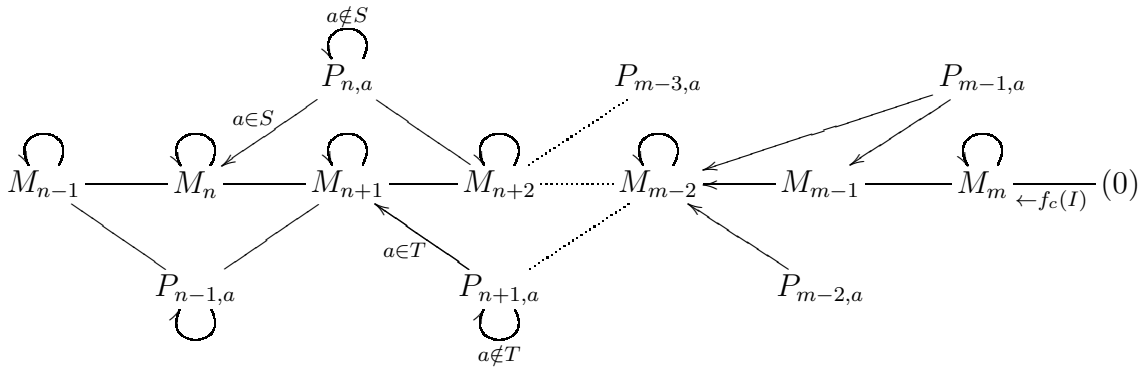
\xymatrix{
                        &          & P_{n,a}\ar @(ru,lu)[]_{a \notin S} \ar @{->}_{a \in S} [dl] \ar @{-} [dr]&          & P_{m-3,a} \ar @{.} [dl]&         & P_{m-1,a} \ar @{->} [dl] \ar @{->} [dll]& \\
     M_{n-1} \ar @(ru,lu)[] \ar @{-} [r] \ar @{- } [dr] & M_{n}\ar @(ru,lu)[]\ar @{-} [r]& M_{n+1}\ar @(ru,lu)[]\ar @{-} [r] \ar @{-} [dr]&M_{n+2}\ar @(ru,lu)[]\ar @{.} [r]&M_{m-2}\ar @(ru,lu)[] \ar @{-} [r]&M_{m-1}\ar @{->} [l] &M_m\ar @(ru,lu)[] \ar @{-} [l]&(0)\ar @{-}^{\leftarrow f_c(I)} [l]\\
               & P_{n-1,a}\ar @(rd,ld)[] \ar @{-} [ur]     &      &  P_{n+1,a}\ar @(rd,ld)[]^{a \notin T} \ar @{->}^{a \in T} [ul] \ar @{.} [ur]  &       &P_{m-2,a} \ar @{->} [ul]     &  & } $$
\qed

Surprisingly, the semiprime operations of the form $f_{n,S,T,m\prime}^f$
and $g_{n,S,T,m\prime}^f$ do not commute with some of the other semiprime
operations for a non-zero ideal. For example,
$$f_{n,S,T,m\prime}^f\circ  f_{n,S,T,m-1\prime}^f(M_m)=M_{m-2} \textrm{ but
} f_{n,S,T,m-1\prime}^f\circ f_{n,S,T,m\prime}^f(M_m)=M_{m-1}.$$
Also,
$$f_{n,S,T,m\prime}^f\circ f_{m-1,a}^f(M_{m+1})=M_{m-2} \textrm{ but
} f_{m-1,a}^f\circ f_{n,S,T,m\prime}^f(M_{m+1})=M_{m-1}.$$  This makes it
hard to decompose the semiprime operations of $R=K[[t^2,t^3]]$, $S_R$,
into the union of submonoids of $M_{\mathfrak{I}}$ like we did in the
Dedekind case.

We make the following definition:

\begin{dfn}  Let $R$ be a one-dimensional semigroup ring defined by
$S \subseteq \mathbb{N}_0$.  Let $f_c$ is a bounded semiprime operation
and $J$ be the unique ideal with $f_c(I)=J$ for all $(0) \neq I  \subseteq
J$ and $n \geq 1$ be the conductor of $S$. Suppose $\mathfrak{a}$ is an
ideal which is incomparable to $J$ and $f_c(\mathfrak{a}) \supseteq J$ and
$\mathfrak{a}=\mathfrak{a}_0 \subseteq \mathfrak{a}_1 \subseteq \cdots
\subseteq \mathfrak{a}_k=f_c(\mathfrak{a})$ is a composition series for
$f_c(\mathfrak{a})/\mathfrak{a}$ for $k \geq n$ with $\mathfrak{a}_i
\supseteq J$ for all $i>0$. Then we say $f_c$ is an exceptional semiprime
operation.
\end{dfn}

Note that the semiprime operations $f_{n,S,T,m\prime}^f$ and
$g_{n,S,T,m\prime}^f$ are exceptional bounded semiprime operations
since $P_{m-1,b} \subseteq M_{m-1} \subseteq M_{m-2}$ is a
composition series for $f_c(P_{m-1,b})/P_{m-1,b}$ of length $2$
which is the conductor of $<2,3>$ the semigroup associated to
$R=K[[t^2,t^3]]$.

Combining the results of Theorems \ref{usp} and \ref{bsp} and looking at
compositions of the maps obtained in the theorems we see that the
non-exceptional semiprime operations can be decomposed as in the Dedekind
case:

\begin{thrm} Let $R=K[[t^2,t^3]]$ and $E$ be the set of exceptional semiprime operations of $R$.
Then the complement of $E$ in $S_R$,
$S_R \setminus E$, is the union of the monoids
$$M_0=\{e\} \cup \{f_{n,S,T}^{int},f_{n,a}^f, f_{n,S,T,m}^f  \}$$ and $$M_f=\{e\}
 \cup \{g_{n,a}^f, g_{n,S,T,m}^f \}$$ where $M_f$ is a left $M_0$-act but not a right
$M_0$-act under composition.
\end{thrm}

{\bf Proof:} Above we saw by example that the semiprime operations
$f_{n,S,T,m\prime}^f$ and $g_{n,S,T,m\prime}^f$ were exceptional. To
see the remaining bounded semiprime operations are not exceptional,
we need to find all nonzero ideals $\mathfrak{a}$ which are not
comparable to the ideal $J$ for each bounded semiprime operation
$f_c$ for which $f_c(I)=J$.

For both $f_c=f_{n,S,T,m}^f$ and $f_c=g_{n,S,T,m}^f$, the $J$ in the
definition is $M_m$.  The only ideals which are incomparable to
$M_m$ are $P_{m-1,a}$ for all $a \in K$ and
$f_c(P_{m-1,a})=M_{m-1}\supseteq M_m$ and $P_{m-1,a} \subseteq
M_{m-1}$ is a composition series for $f_c(P_{m-1,a})/P_{m-1,a}$.
Thus $f_{n,S,T,m}^f$ and $g_{n,S,T,m}^f$ are not exceptional.

For both $f_c=f_{m,a}^f$ and $f_c=g_{m,a}^f$, the $J$ in the
definition is $P_{m,a}$.  The ideals which are incomparable to
$P_{m,a}$ are $P_{m,b}$ for $b \neq a$, $M_{m+1}$, $P_{m+1,b}$ for
all $b \in K$ and $P_{m-1,b}$ for all $b \in K$.  Note that
$f_c(P_{m,b})=M_m \supseteq P_{m,a}$ and $P_{m,b} \subseteq M_m$ is
a composition series for $f_c(P_{m,b})/P_{m,b}$.  Also
$f_c(M_{m+1})=M_m \supseteq P_{m,a}$ and $M_{m+1} \subseteq M_{m}$
is a composition series for $f_c(M_{m+1})/M_{m+1}$.  Similarly,
$f_c(P_{m+1,b})=M_m \supseteq P_{m,a}$ and $P_{m+1,b} \subseteq
M_{m+1} \subseteq M_{m}$ is a composition series for
$f_c(P_{m+1,b})/P_{m+1,b}$ and only $M_m$ is comparable to
$P_{m,a}$.  Lastly, $f_c(P_{m-1,b})=M_{m-1} \supseteq P_{m,a}$ and
$P_{m-1,b} \subseteq M_{m-1}$ is a composition series for
$f_c(P_{m-1,b})/P_{m-1,b}$.  Now by definition both $f_{m,a}^f$ and
$g_{m,a}^f$ are not exceptional.

Now we look at all compositions of semiprime operations in $M_0$.
Throughout, we will denote $K \backslash \{a\}=a^{C}$. The
compositions are as follows:
\begin{enumerate}[(M1)]

\item  $f_{m,S,T}^{int} \circ f_{n,U,V}^{int}=f_{n,U,V}^{int} \circ
f_{m,S,T}^{int}=\left\{ \begin{array}{l}f_{m,S,T}^{int} \textrm{ if } m+2 \leq n\\
                             f_{m,S,T \cup U}^{int} \textrm{ if } m+1 =n\\
                             f_{m,S \cup U,T \cup V}^{int} \textrm{ if } m= n\\
                             f_{n,U,V \cup S}^{int} \textrm{ if } n+1 =m\\
                             f_{n,U,V}^{int} \textrm{ if } n+2 \leq m \end{array}\right.$

\item $f_{n,S,T,m}^f \circ f_{l,U,V}^{int} =f_{l,U,V}^{int}\circ f_{n,S,T,m}^f
             =\left\{ \begin{array}{l} f_{n,S,T,m}^f \textrm{ if } n+1 \leq m < l\\
                            f_{n,S,T \cup U,m}^f \textrm{ if } n+1=l \leq m\\
                            f_{n,S \cup U,T \cup V,m}^f \textrm{ if } n=l \leq m-1\\
                            f_{l,U,S \cup V,m}^f \textrm{ if } l+1=n \leq m-1\\
                            f_{l,U,V ,m}^f \textrm{ if } l+1 <n \leq m-1\\
                       \end{array}\right.$

\item $f_{m,a}^f \circ f_{l,U,V}^{int} =f_{l,U,V}^{int}\circ
f_{m,a}^f=\left\{ \begin{array}{l} f_{m,a}^f \textrm{ if } m < l, m=l, a \notin U \textrm{ or } l=m-1, a \notin V\\
                        f_{{m-1},K,K,m}^f \textrm{ if } m=l, a \in U\\
                       f_{m-1,U,K,m}^f \textrm{ if } l= m-1, a \in V\\
                        f_{l,U,V,m}^f \textrm{ if } l < m-1\\
                       \end{array}\right.$

\item $f_{n,S,T,m}^f \circ f_{l,U,V,k}^{f} =f_{l,U,V,k}^{f}\circ
f_{n,S,T,m}^f=\left\{ \begin{array}{l} f_{n,S,T,m}^f \textrm{ if } n+1 < m,l\\
                       f_{n,S,T \cup U,m}^f \textrm{ if } n+1=l\leq m<k\\
                       f_{n,S,T \cup U,k}^f \textrm{ if } n+1=l \leq k-1\leq m-1\\
                       f_{n,S\cup U,T\cup V,m}^f \textrm{ if } n=l <m-1 <k-1\\
                       f_{n,S\cup U,T\cup V,k}^f \textrm{ if } n=l <k-1 \leq m-1\\
                       f_{l, U,S\cup V,m}^f \textrm{ if } l+1=n \leq m-1 <k-1\\
                       f_{n,U,S\cup V,k}^f \textrm{ if } l+1=n <k \leq m\\
                       f_{l,U,V,k}^f \textrm{ if } l+1<m,k\\
                       \end{array}\right.$

\item $f_{n,S,T,m}^f \circ f_{l,a}^{f} =f_{l,a}^{f}\circ
f_{n,S,T,m}^f= \left\{ \begin{array}{l} f_{n,S,T,m}^f \textrm{ if } m \leq l \\
                       f_{n,S,T,l}^f \textrm{ if } n+1<l\leq m\\
                       f_{n,K,T \cup {a}^C,l}^f \textrm{ if } n+1=l \leq m\\
                       f_{n-1,K,K,n}^f \textrm{ if } n=l \leq m-1\\
                       f_{l,a}^f \textrm{ if } l <n\\
                       \end{array}\right.$

\item $f_{n,a}^f \circ f_{m,b}^f= f_{m,b}^f \circ f_{n,a}^f=
\left\{ \begin{array}{l}f_{n,a}^f \textrm{ if } n+1<m \\
                       f_{n-1,K,K,n}^f \textrm{ if } m \leq n \leq m+1\\
                       f_{m-1,K,K,m}^f \textrm{ if } n+1=m\\
                       f_{m,b}^f \textrm{ if } m+1<n\\
                       \end{array}\right.$

\end{enumerate}
Clearly, $M_0$  is a monoid and similar compositions show that $M_f$ is a
monoid.  To see that $M_f$ is a left $M_0$-act but not a right $M_0$-act
we look at the mixed compositions.

\begin{enumerate}[(L1)]
\item  $g_{n,S,T,m}^f \circ f_{l,U,V}^{int} =f_{l,U,V}^{int}\circ
g_{n,S,T,m}^f=\left\{ \begin{array}{l} g_{n,S,T,m}^f \textrm{ if } n+1 \leq m, l\\
                       g_{n,S,T \cup U,m}^f \textrm{ if } n+1=l \leq m\\
                       g_{n,S\cup U,T \cup V,m}^f \textrm{ if } n=l \leq m-1\\
                       g_{l,U,S \cup V,m}^f \textrm{ if } l+1=n \leq m-1\\
                       g_{l,U,V,m}^f \textrm{ if } l < n \leq m-1\\
                       \end{array}\right.$

\item $g_{m,a}^f \circ f_{l,U,V}^{int} =f_{l,U,V}^{int}\circ
g_{m,a}^f=\left\{ \begin{array}{l} g_{m,a}^f \textrm{ if } m < l, m=l, a \notin U \textrm{ or } l=m-1, a \notin V\\
                        g_{{m-1},K,K,m}^f \textrm{ if } m=l, a \in U\\
                       g_{m-1,U,K,m}^f \textrm{ if } l= m-1, a \in V\\
                        g_{l,U,V,m}^f \textrm{ if } l < m-1\\
                       \end{array}\right.$

\item

(a) $g_{n,S,T,m}^f \circ f_{l,U,V,k}^{f} =
            \left\{ \begin{array}{l} g_{n,S,T,m}^f \textrm{ if } n+1 < l, m \leq k\\
                          g_{n,S,T \cup U,m}^f \textrm{ if } n+1=l < m \leq k  \\
                          g_{n,S\cup U,T \cup V,m}^f \textrm{ if } n=l \leq m-1 \leq k-1  \\
                          \textrm{ not a semiprime operation if }  k < m \\
                          \end{array}\right.$

\noindent(b) $f_{l,U,V,k}^{f} \circ g_{n,S,T,m}^f  =\left\{ \begin{array}{l} g_{n,S,T,m}^f \textrm{ if } n+1 < l, m \leq k\\
                          g_{n,S,T \cup U,m}^f \textrm{ if } n+1=l < m \leq k  \\
                          g_{n,S,T \cup U,k}^f \textrm{ if } n+1=l \leq k < m  \\
                          g_{n,S\cup U,T \cup V,m}^f \textrm{ if } n=l \leq m-1 \leq k-1  \\
                          g_{n,S\cup U,T \cup V,k}^f \textrm{ if } n=l < k \leq m  \\
                          g_{l,U,V \cup S,m}^f \textrm{ if } l+1=n < m \leq k  \\
                          g_{l,U,V \cup S,k}^f \textrm{ if } l+1=n \leq k < m  \\
                          g_{l,U,V,k}^f \textrm{ if } l+1 < n, k \leq m\\
                          \end{array}\right.$

\item
(a) $g_{n,a}^f \circ f_{m,b}^f=\left\{ \begin{array}{l}g_{n,a}^f \textrm{ if } n+1<m \\
                       \textrm{ not a semiprime operation if }  m \leq n+1 \\
                       \end{array}\right.$

\noindent(b)  $f_{m,b}^f \circ g_{n,a}^f= \left\{ \begin{array}{l}g_{n,a}^f \textrm{ if } n+1<m \\
                       g_{n-1,K,K,n}^f \textrm{ if } m \leq n \leq m+1\\
                       g_{m-1,K,K,m}^f \textrm{ if } n+1=m\\
                       g_{m,b}^f \textrm{ if } m+1<n\\
                       \end{array}\right.$

\item (a) $g_{n,S,T,m}^f \circ f_{l,a}^{f} = \left\{ \begin{array}{l} g_{n,S,T,m}^f \textrm{ if } m \leq l \\
                       \textrm{ not a semiprime operation if }  l < m \\
                       \end{array}\right.$

\noindent(b) $  f_{l,a}^{f}\circ g_{n,S,T,m}^f = \left\{ \begin{array}{l} g_{n,S,T,m}^f \textrm{ if } m \leq l \\
                       g_{n,S,T,l}^f \textrm{ if } n+1<l\leq m\\
                       g_{n,K,T \cup {a}^C,l}^f \textrm{ if } n+1=l \leq m\\
                       g_{n-1,K,K,n}^f \textrm{ if } n=l \leq m-1\\
                       g_{l,a}^f \textrm{ if } l <n\\
                       \end{array}\right.$

\item (a) $  g_{l,a}^{f}\circ f_{n,S,T,m}^f = \left\{ \begin{array}{l}
                       g_{l,a}^f \textrm{ if } l <n\\
                       \textrm{ not a semiprime operation if }  l \geq n \\
                       \end{array}\right.$

\noindent(b) $f_{n,S,T,m}^f \circ g_{l,a}^{f} = \left\{ \begin{array}{l} g_{n,S,T,m}^f \textrm{ if } m \leq l \\
                       g_{n,S,T,l}^f \textrm{ if } n+1<l\leq m\\
                       g_{n,K,T \cup {a}^C,l}^f \textrm{ if } n+1=l \leq m\\
                       g_{n-1,K,K,n}^f \textrm{ if } n=l \leq m-1\\
                       g_{l,a}^f \textrm{ if } l <n\\
                       \end{array}\right.$

\end{enumerate}
Hence $M_f$ is a left $M_0$-act but not a right $M_0$-act.  \qed

We will now see as in the Dedekind case the only prime operation is the
identity.

\begin{thrm} Let $R=K[[t^2,t^3]]$.  Then $P_R=\{e\}$.
\end{thrm}

{\bf Proof:}  Suppose $f_c$ is one of the other semiprime
operations. Then for some $i\geq 2$ and some $a \in K$,
$f_c(P_{i,a})=M_{i}$. Now since $f_c$ is prime,
$P_{i,a}=(t^i+at^{i+1})f_c(R)=f_c(P_{i,a})=M_i$ which is a
contradiction. Hence $f_c$ cannot be prime.  Thus $P_R=\{e\}$. \qed

To determine all the semiprime operations for other semigroup rings
becomes immediately more complicated for any other semigroup. Note even
for the ring $K[[t^2,t^5]]$, the diagram of two generated monomial ideals
is as follows:

$$ 
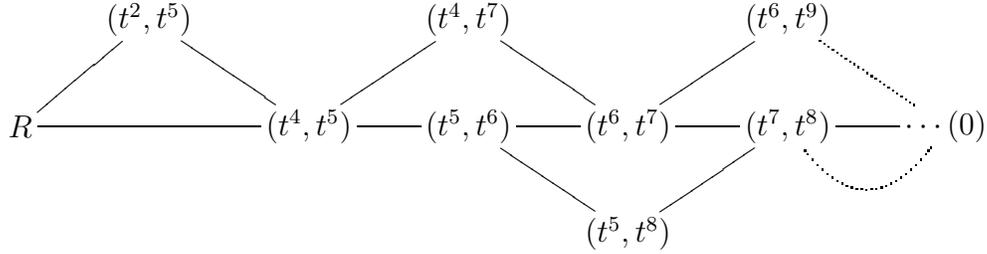
\xymatrix{
                 & (t^2,t^5)\ar @{-} [dr]      &           & (t^4,t^7) \ar @{-} [dl] \ar @{-} [dr]&          & (t^6,t^{9}) \ar @{-} [dl]\ar @{.} [dr]&         & \\
    R\ar @{-} [rr] \ar @{-} [ur]& & (t^4,t^5)\ar @{-} [r]& (t^5,t^6)\ar @{-} [r] \ar @{-} [dr]&(t^6,t^7)\ar @{-} [r]&(t^7,t^8) \ar @{-} [r]&\cdots  (0)\ar @/^2pc/ @{.} [l]\\
                 &       &    &      &   (t^5,t^8) \ar @{-} [ur]     &       &       } $$
where each line segment in the above diagram indicates $\supseteq$.
Of course, this leaves out a lot of two generated ideals in addition
to all the principal ideals.  But even without all these ideals we
can see that there is an extra layer of difficulty that we did not
have in the cuspidal cubic case.  Certainly, the conductor will be
involved with the classification of all semiprime operations.  I
believe that the non-exceptional semiprime operations over a
one-dimensional semigroup ring $R$ will decompose into the union of
two submonoids of the monoid $(M_{\mathfrak{I}},\circ)$ of maps from
the set of ideals of $R$ to itself, one being a left but not a right
act of the other.

Certainly, if $f_c$ is a prime operation over any commutative ring,
then $f_c$ is the identity on the set of principal ideals of $R$
since $gR=gf_c(R)=f_c(g)$  for all $g \in R$.  However, it is not
known whether $f_c$ must be the identity over one dimensional
domains. It may be that for one-dimensional semigroup rings, the set
of prime operations will be the singleton set consisting only of the
identity.

There will certainly be more prime operations if the ring is a
normal domain of dimension 2 or more since the integral closure does
not agree with the identity for all ideals of height 2 or more.
Moreover, integral closure is a prime operation in any normal
domain.

\end{document}